\documentclass[a4paper,11pt]{article}

\usepackage{amsfonts, amsmath, amssymb, amsthm}
\usepackage[colorlinks=true,citecolor=blue]{hyperref}
\usepackage[margin = 2.50cm]{geometry}
\usepackage{mathtools}
\usepackage{graphicx}
\usepackage{enumerate}
\usepackage{verbatim}
\usepackage{tikz}
\usepackage{algorithmic}
\usepackage[utf8]{inputenc}
\usepackage{microtype}
\usepackage{amsmath,amsthm,amssymb,color,tikz,mathtools}

\usepackage{tcolorbox}
\usepackage{mdframed}

\usetikzlibrary{positioning}
\usetikzlibrary{shapes}
\usetikzlibrary{decorations.pathreplacing}
\usepackage{enumerate, enumitem}

\newcommand{\FF}{\mathbb{F}}

\newcommand{\NN}{\mathbb{N}}

\newcommand{\RR}{\mathbb{R}}

\newcommand{\ZZ}{\mathbb{Z}}

\newcommand{\cQ}{\mathcal{Q}}

\newcommand{\set}[1]{\left\{ #1 \right\}}

\newcommand{\zone}{\set{0,1}}

\renewcommand{\epsilon}{\varepsilon}

\renewcommand{\deg}{\text{deg}}

\newtheorem{theorem}{Theorem}
\newtheorem{corollary}[theorem]{Corollary}
\newtheorem{remark}[theorem]{Remark}

\newtheorem{defi}[theorem]{Definition}

\newtheorem{proposition}[theorem]{Proposition}

\newtheorem*{theorem*}{Theorem}

\newcommand{\remove}[1]{}
\newcommand{\ic}{\mathrm{ic}}
\newcommand{\w}{\mathrm{wt}}
\newcommand{\spc}{\mathrm{spc}}

\begin{document}

\title{
About almost covering subsets of the hypercube
}

\author{
Arijit Ghosh 
\footnote{Indian Statistical Institute, Kolkata, India}
\and
Chandrima Kayal 
\footnote{Universit\'e Paris Cit\'e, CNRS, IRIF, Paris, France}
\and 
Soumi Nandi 
\footnote{The Institute of Mathematical Sciences, Chennai, India}
}

\date{}

\maketitle

\begin{abstract}
Let $\FF$ be a field and consider the hypercube $\{0,1\}^{n}$ in $\FF^{n}$. Sziklai and Weiner (Journal of Combinatorial Theory, Series A 2022) showed that if a polynomial $P(X_{1}, \dots, X_{n}) \in \FF [X_{1}, \dots, X_{n}]$ vanishes on every point of the hypercube $\{0,1\}^{n}$ except those with at most $r$ many ones then the degree of the polynomial will be at least $n-r$. This is a generalization of Alon and F\"uredi's fundamental result (European Journal of Combinatorics 1993) about polynomials vanishing on every point of the hypercube except at the origin (point with all zero coordinates). 
Sziklai and Weiner proved their result using the M\"{o}bius inversion formula and the Zeilberger method for proving binomial equalities. In this short note, we show that a stronger version of Sziklai and Weiner's result can be derived directly from Alon and F\"uredi's result. We also prove a multiplicity version of our results when $\FF = \RR$.
\end{abstract}

\section{Introduction}

A classical result of Alon and F\"{u}redi, also derivable from Alon’s Combinatorial Nullstellensatz~\cite{Alon99}, gives a sharp degree lower bound for polynomials vanishing on almost all points of the $n$-dimensional hypercube $\zone^n$ embedded in $\FF^n$, where $\FF$ is any field.
     
\begin{theorem}[Alon and F\"{u}redi~\cite{AF93}]
   Suppose $P(X_{1}, \dots, X_{n})$ is a polynomial in $\FF \left[ X_{1}, \dots, X_{n} \right]$ such that  $P(u)=0$ for all $u\in\zone^n\setminus\{(0,\dots,0)\}$ and $P(0,\dots,0)\neq 0$. Then $deg(P) \geq n$. 
   \label{th: AF}
\end{theorem}

This result has led to several generalizations and applications. A natural question is: for a proper subset $V \subseteq \zone^n$, what is the minimum degree of a polynomial $P \in \FF[X_1,\dots,X_n]$ that vanishes on $V$ but not on the entire cube?  
For $u=(u_1,\dots,u_n)\in\zone^n$, let the \emph{weight} $\w(u)$ denote the number of coordinates equal to~1.  
Using the M\"{o}bius inversion formula~\cite[Chapter~3]{Stanley_2011} together with Zeilberger’s algorithm for Binomial identities~\cite{PauleS95}, Sziklai and Weiner obtained the following extension of Theorem~\ref{th: AF}.

\begin{theorem}[Sziklai and Weiner~\cite{SziklaiW23}]
Let $R(X_{1},\dots,X_{n}) \in \FF[X_1,\dots,X_n]$ be such that $R(u)\neq 0$ for all $u\in\zone^n$ with $\w(u)\leq r$, and $R(v)=0$ for all $v\in\zone^n$ with $\w(v)>r$. Then $\deg(R)\geq n-r$.
\label{eqvSW23}
\end{theorem}

This can also be expressed in the following equivalent form.

\begin{theorem}[Reformulation of Theorem~\ref{eqvSW23}]
Suppose $P(X_{1},\dots,X_{n}) \in \FF[X_1,\dots,X_n]$ satisfies $P(u)=0$ for all $u\in\zone^n$ with $\w(u)\leq k$, and $P(v)\neq 0$ for all $v\in\zone^n$ with $\w(v)>k$. Then $\deg(P)>k$.
\label{SW23}
\end{theorem}

The equivalence follows by considering
\[
R(X_1,\dots,X_n) = P(1-X_1,\dots,1-X_n),
\]
which transforms small-weight conditions into large-weight ones without changing the degree.  
Using Gr\"{o}bner basis methods, Heged\"{u}s further strengthened this result.

\begin{theorem}[Heged\"{u}s~\cite{hegedus2024coversfinitesetspoints}]
Suppose $P(X_{1},\dots,X_{n}) \in \FF[X_1,\dots,X_n]$ satisfies $P(u)=0$ for all $u\in\zone^n$ with $\w(u)\leq k$, and there exists $v\in\zone^n$ with $\w(v)>k$ such that $P(v)\neq 0$. Then $\deg(P)>k$.
\label{thm-Hegedus}
\end{theorem}

\paragraph{Our contribution.}
In this short note, we show that a stronger result, encompassing Theorems~\ref{SW23} and~\ref{thm-Hegedus}, follows directly from Alon and Füredi’s theorem.


\begin{theorem}[Main result]
Let $P \in \FF[X_1,\dots,X_n]$ be a polynomial. Suppose there exist $a,b \in \zone^n$ such that $P(a)=0$ and $P(b)\neq 0$. Then
\[
\deg(P) \;\ge\; \max\{w,\, n-W\},
\]
where
$w = \min \left\{\w(z) : z \in \zone^n,\; P(z)\neq 0 \right\}$ and $W = \max \left\{\w(z) : z \in \zone^n,\; P(z)\neq 0 \right\}$.
\label{thm-new-lower-bound-new}
\end{theorem}

\color{black}

Theorem~\ref{thm-Hegedus} follows as a direct corollary of the above result, and since Theorem~\ref{thm-Hegedus} generalizes Theorem~\ref{SW23}, our result subsumes both.  
Setting $w=W=k$ gives the following.

\begin{corollary}
Suppose $P\in\FF[X_1,\dots,X_n]$ satisfies $P(u)\neq 0$ if and only if $\w(u)=k$, for some $k\in\{0\}\cup [n]$. Then $\deg(P)\geq \max\{k,\,n-k\}$.
\label{cor1}
\end{corollary}

For $\FF=\RR$, this bound was proved independently by Venkitesh (Corollary 6.2~\cite{Venkitesh22}) and by Ghosh, Kayal, and Nandi (Theorem 1.3~\cite{GKN23}), with the latter also showing tightness.

\medskip

A further consequence of Theorem~\ref{thm-new-lower-bound-new} concerns Boolean functions over $\FF_2^n$, that is, maps from $\FF_2^n$ to $\FF_2$, whose algebraic degree plays a key role in cryptography, pseudorandomness, and coding theory.

\begin{corollary}\label{cor:degree_Boolean_functions}
Suppose $P\in\FF_2[X_1,\dots,X_n]$ vanishes at $v\in\zone^n$ if and only if $\w(v)=k$. Then $\deg(P)\geq \max\{k,\,n-k\}$.
\end{corollary}

This bound holds only over $\FF_2$. For other prime fields it can fail, as the following example shows. When $p \geq 3$ and $n=p$, the polynomial $P(X_1,\dots,X_n)=\sum_{i=1}^n X_i - 2$ over $\FF_p$ vanishes exactly on points of weight $2$ on $\zone^n$, yet has degree $1<\max\{2,p-2\}$.




\paragraph{Multiplicity version.}
A polynomial $P \in \RR[X_1,\dots,X_n]$ has a \emph{zero of multiplicity $\ell$} at $a \in \RR^n$ if $P(a)=0$ and all partial derivatives of $P$ up to order $\ell-1$ vanish at $a$.
We consider a natural extension: for a subset $S \subsetneq \zone^n$ and an integer $\ell \geq 2$, what is the minimum degree of a real polynomial that vanishes on every point of $S$ with multiplicity at least $\ell$ and is nonzero on some point of $\zone^n \setminus S$?



Clifton and Huang~\cite{CH20} introduced hyperplane coverings with multiplicities of the hypercube excluding the origin. Sauermann and Wigderson~\cite{SW20} subsequently initiated the study of polynomial coverings in this higher-multiplicity setting. Ghosh, Kayal, Nandi, and Venkitesh~\cite{GKNV23} further investigated polynomial coverings of \emph{symmetric}\footnote{A subset $S \subseteq \{0,1\}^n$ is symmetric if, whenever it contains a point of weight $k$, it contains all points of weight $k$ in $\{0,1\}^n$.} subsets of the hypercube with multiplicities.

\color{black}

\begin{theorem}
\label{th:mult_R_new}
Let $\ell \ge 2$, and $S \subsetneq \zone^{n}$ be a non-empty subset of $\zone^{n}$. Suppose $P\in\RR[X_1,\dots,X_n]$ is a nonzero polynomial that vanishes on $S$ with multiplicity at least $\ell$. Assume there exist integers $0 \le w_\ell \le W_\ell \le n$ such that:
\begin{itemize}
    \item there exists $u \in \zone^{n}$ with $\w(u)=w_\ell$ and $P(u)\ne 0$;
    \item there exists $v \in \zone^{n}$ with $\w(v)=W_\ell$ and $P(v)\ne 0$; and
    \item for all $y \in \zone^{n}$, if $\w(y) < w_\ell$ or $\w(y) > W_\ell$, then $P$ has a zero of multiplicity at least $\ell$ at $y$.
\end{itemize}
If 
\[
\max\left\{ w_\ell,\, n-W_\ell \right\} \geq 2\ell-3,
\]
then
\[
\deg (P) \geq \max\left\{ w_\ell,\,n-W_\ell \right\}+2\ell-3.
\]
\label{th: mult_R}
\end{theorem}
\noindent
In Proposition~\ref{proposition:lower_bound_tightness} we establish the tightness of the above lower bound.

\color{black}

As a direct consequence of Theorem~\ref{th:mult_R_new} we obtain the following extension of the non-multiplicity lower bound for the {\em $k$-th layer} of the hypercube proved in~\cite{Venkitesh22,GKN23}.

\begin{corollary}
Let $\ell\in \NN$ and $k\in \{0\}\cup [n]$. Suppose $P\in\RR[X_1,\dots,X_n]$ is nonzero, vanishes with multiplicity at least $\ell$ on every $v\in\zone^n$ with $\w(v)\neq k$, and is nonzero on every $u\in\zone^n$ with $\w(u)=k$. If
$\max\left\{k,\,n-k\right\} \geq 2\ell-3$,
then
\[
\deg(P) \geq \max \left\{k,\,n-k \right\}+2\ell-3.
\]
\end{corollary}

\noindent
For $k=0$, Sauermann and Wigderson~\cite{SW20} constructed a polynomial of degree $n+2\ell-3$ that is nonzero at the origin and vanishes with multiplicity at least $\ell$ at every other point of the hypercube, showing that the bound above is tight in this case. For $k\geq 1$, the existence of matching polynomials remains a compelling open problem.

\paragraph{Definitions and notation.}
We will use the following definitions and notation in this paper. 
\begin{itemize}

    \item 
        We will denote by $\RR$ the set of all real numbers, $\NN$ the set of all natural numbers, and $\mathbb{Z}$ the set of all integers.

    \item 
        Given a set $S$, $\binom{S}{k}$ denotes the set of all $k$-sized subsets of $S$.
    
    \item 
        For any $n \in \NN$, $[n]$ denotes the set $\left\{ 1, \dots, n\right\}$.

    \item 
        Given $n \in \NN$, $\mathbb{Z}_{n}$ denotes the ring $\mathbb{Z}/{\left( n\mathbb{Z} \right)}$, and $\FF_{2}$ denotes the field with two elements. 

    \item 
        $\mathcal{R}[X_{1}, \dots, X_{n}]$ denotes the polynomial ring over the ring $\mathcal{R}$. 


    \item 
        Given a polynomial $P(X_{1}, \dots, X_{n})$ over the polynomial ring $\mathcal{R}[X_{1}, \dots, X_{n}]$, $deg(P)$ denotes the degree of the polynomial $P(X_{1}, \dots, X_{n})$.
    
    \item 
        Given $u \in \{ 0,1\}^{n}$, $\w(u)$, {\em weight} of $u$, denotes the $L_{1}$-norm $\|u\|_{1}$ of $u$, that is, number of $1$'s in the coordinates of $u$.

\end{itemize}

\section{Covering with polynomials from $\FF[X_1, \dots, X_n]$}

We now proceed to prove Theorem~\ref{thm-new-lower-bound-new}.
\begin{proof}[Proof of Theorem~\ref{thm-new-lower-bound-new}]
Note that $P$ is non-constant and hence $\deg(P) > 0$. If $w = 0$ and $W = n$, then $\max\{w,\, n-W\} = 0$, and the statement holds trivially. Thus, we may assume that this is not the case; in particular, we can take
\[
    \max\{w,\, n-W\} \ge 1.
\]
We may further assume that $w \ge n - W$; otherwise, we replace $P$ by the polynomial
\[
\widetilde{P}(X_1,\dots,X_n) := P(1 - X_1,\dots,1 - X_n).
\]
Therefore, $w\ge 1$.
Let $u = (u_{1}, \dots, u_{n}) \in \zone^{n} \setminus S$ be such that $\w(u) = w$ and $P(u)\neq 0$. We know that, for all $v\in\zone^n$ with $\w(v) < w$, we have $P(v)=0$. Without loss of generality we may assume that $u_i=1$ if and only if $i\in[w]$. Consider the polynomial 
$$
    Q(X_1,\dots,X_w) : = P(1-X_1,\dots, 1-X_w,0,\dots,0). 
$$
Then $\deg(Q)\leq\deg(P)$ and $Q(0,\dots,0)=P(u)\neq 0$. 
Now take any $y\in\zone^w\setminus\{(0,\dots,0)\}$. Then $Q(y)=P(\Tilde{y})$, where $\Tilde{y}_{i}=1-y_i$ for all $i\in [w]$ and $\Tilde{y}_{i}=0$ for all $i\in [n]\setminus [w]$. Since $\w(\Tilde{y})<w$, we have $P(\Tilde{y})=0$. So the polynomial $Q(X_{1}, \dots, X_{w})$ vanishes everywhere on $\zone^w$, except at $(0,\dots,0)\in\zone^w$. Therefore by Theorem~\ref{th: AF}, $\deg(Q)\geq w$. Hence, $\deg(P)\geq w$.
\end{proof}

\color{black}
We will now give a proof of Corollary~\ref{cor:degree_Boolean_functions}.

\begin{proof}[Proof of Corollary~\ref{cor:degree_Boolean_functions}]
    Let $Q\in\FF_2[X_1,\dots,X_n]$ be defined by $Q(X_1,\dots,X_n)=1-P(X_1,\dots,X_n)$. Then $\deg(Q)=\deg(P)$ and for any $u\in\zone^n$, $Q(u)=0$ if and only if $\w(u)\neq k$.
    Then by Theorem~\ref{thm-new-lower-bound-new}, $\deg(Q)\geq\max\{k,n-k\}$.  
\end{proof}
The next proposition establishes the tightness of Corollary~\ref{cor1}, and consequently of Corollary~\ref{cor:degree_Boolean_functions}.

\begin{proposition}
    Let $q\in \NN$ be a prime, and $0\leq k \leq n$ satisfy the following: there exists $m\in\NN$ such that 
    \begin{itemize}
        \item $q^{m-1}\leq n-k<q^m$ and
        \item $k+1$ is divisible by $q^m$.
    \end{itemize}
    Then there exists a polynomial $P(X_{1}, \dots, X_{n}) \in\FF_q[X_1,\dots,X_n]$ such that $\deg(P)=k$ and for any $u\in\zone^n$, $P(u)=0$ if and only if $\w(u)\neq k$. 
\label{prop: tightness}\end{proposition}

Note that, as $k+1$ is divisible by $q^m$, $k+1\geq q^m> n-k$, which implies $k\geq\frac{n}{2}$, that is, $k=\max\{k,n-k\}$.
We need the following result by Kummer~\cite{Kummer1852} to prove Proposition~\ref{prop: tightness}.

\begin{theorem}[Kummer~\cite{Kummer1852}]
    Let $j$ and $m$ be integers such that $0\leq j\leq m$. Let $\alpha\in\NN$ and $p$ be a prime. Then $p^\alpha$ divides $\binom{m}{j}$ if and only if at least $\alpha$ carries are needed when adding $j$ and $m-j$ in base $p$.
\label{kummer}\end{theorem}

\begin{proof}[Proof of Proposition~\ref{prop: tightness}]
    
    We shall show that there exists $P_{k}\in\FF_q[X_1,\dots,X_n]$ with $\deg(P_k)=k$ such that for any $v\in\zone^n$ $P(v)\neq 0$ if and only if $\w(u) = k$. We define 
    $$
        P_{k}(X_1,\dots, X_n) := \sum_{I\in\binom{[n]}{k}}\left(\prod_{i\in I}X_i \right).
    $$
    Then $P_{k}\in\FF_q[X_1,\dots,X_n]$ satisfies the following:
    \begin{itemize}
        \item 
            $\deg(P_{k})=k$, 
        
        \item 
            $P_{k}(\Tilde{u})=0$, for all $\Tilde{u}\in\zone^n$ with $\w( \Tilde{u} )< k$, 
        
        \item 
            $P_{k}(u)=1$, for all $u\in\zone^n$ with $\w( u ) = k$ and 

        \item 
            For all $u_r\in\zone^n$ with $\w( u_{r} )\ = k+r$, where $1\leq r\leq n-k$, we have
            $$
                P_{k}(u_r) \equiv \binom{k+r}{k}  \mod q.
            $$

%
    \end{itemize}

    We claim that for all $r\in [n-k]$, we have
    $$
        \binom{k+r}{k} \equiv 0 \mod q.
    $$
    By our assumption, $k+1$ is divisible by $q^m$. So in the expansion of $k+1$ with powers of $q$, the coefficient of $q^i$ is $0$, for all $i<m$. So in the expansion of $k$ with powers of $q$, coefficient of $q^i$ is $q-1$, for all $i<m$. Now for each $r\in [n-k]\subseteq [q^m-1]$, there exists $j\in\{0,1,\dots,m-1\}$ such that the coefficient of $q^j$ in the expansion of $r$ with powers of $q$ is non-zero. So at least one carry is needed when adding $k$ with $r$ in base $q$. So by Theorem~\ref{kummer}, $q$ divides $\binom{k+r}{k}$, as required.
\end{proof}  



Note that in circuit complexity, one often studies polynomials in $\mathbb{Z}_N[X_1,\dots,X_n]$ that vanish on subsets of the Boolean cube in order to derive degree lower bounds for Boolean functions such as \textsc{OR}~\cite{BBR94} and \textsc{Majority}~\cite{Szegedy89}. Barrington, Beigel, and Rudich~\cite{BBR94} showed that if $P \in \mathbb{Z}_{m}[X_{1}, \dots, X_{n}]$ vanishes only at the origin in $\{0,1\}^{n}$, then $\deg(P) = \Theta\!\bigl(n^{1/r}\bigr)$, where $r$ denotes the number of distinct prime factors of $m$. In contrast, when $m$ is prime, Proposition~\ref{prop: tightness} establishes the existence of a polynomial of a degree $k$ (under suitable conditions) that vanishes at all points of $\{0,1\}^{n}$ except those in the $k$-th layer.

\color{black}
\section{Covering with multiplicities when $\mathbb{F} = \mathbb{R}$}

 In this section we prove Theorem~\ref{th: mult_R}. To prove this result, we will use the following important generalization of Alon and Furedi's result (Theorem~\ref{th: AF}) by Sauermann and Wigderson. 
\begin{theorem}
[Sauermann and Wigderson~\cite{SW20}]\label{prop: Yuval, t,t-1}
Let $\ell\geq 2$, $n\geq 2\ell-3$ and $P(X_1,\dots,X_n)\in\RR[x_1,\dots,x_n]$ be a polynomial  having zeros of multiplicity at least $\ell$ at all points in $\zone^n\setminus\{(0,\dots,0)\}$ and $P(0, \dots, 0)\neq 0$. Then $\deg(P)\geq n+2\ell-3$ and the bound is tight. 
\end{theorem}

\begin{proof}[Proof of Theorem~\ref{th: mult_R}]
    Without loss of generality we assume that $w_{\ell}\geq n-W_{\ell}$, otherwise, we work with the polynomial $\widetilde{P}(X_{1}, \dots, X_{n}) := P(1-X_{1}, \dots, 1-X_{n})$ in place of the polynomial $P$.
Recall that there exists $u \in \zone^{n}$ with $\w(u)=w_\ell$ and $P(u)\ne 0$.
Without loss of generality we may assume that $u_i=1$ if and only if $i\in[w_{\ell}]$. Consider the polynomial $$Q(X_1,\dots,X_{w_{\ell}}) : = P(1-X_1,\dots, 1-X_{w_{\ell}},0,\dots,0).$$ Then $\deg(Q)\leq\deg(P)$ and $Q(0,\dots,0)=P(u)\neq 0$.

Now we claim that $Q(X_{1}, \dots, X_{w_{\ell}})$ vanishes everywhere on $\zone^{w_{\ell}}\setminus\{(0,\dots,0)\}$ with multiplicity at least $\ell$. Take any $a\in\zone^{w_{\ell}}\setminus\{(0,\dots,0)\}$ and let $\Tilde{a}\in\zone^n$ be such that $\Tilde{a}_i=1-a_i$, for all $i\in [w_{\ell}]$ and $\Tilde{a}_i=0$, for all $i>w_\ell$.
Then $\w( \Tilde{a} ) <w_{\ell}$ and so $\Tilde{a}$ is a zero of $P$ with multiplicity at least $\ell$. Hence $Q(a)=P(\Tilde{a})=0$.
Next consider any differential operator $$D=\frac{\partial^{d_1}}{\partial X_1^{d_1}}\dots\frac{\partial^{d_{w_{\ell}}}}{\partial X_{w_{\ell}}^{d_{w_{\ell}}}}$$ such that $d_i\geq 0$, for all $i\in [w_{\ell}]$ and $\sum_{i=1}^{w_{\ell}} d_i< \ell$. We shall show that $DQ(a)=0$.

First observe that we can write $P(X_1,\dots,X_n)=P_1(X_1,\dots,X_{w_{\ell}})+P_2(X_1,\dots,X_n)$, where each monomial of $P_2$ involves at least one $X_i$, where $i\in\{w_{\ell}+1,\dots,n\}$, with degree of $X_i$ at least $1$. So $P_2$ vanishes at each \remove{$(X_1,\dots,X_w,\underbrace{0,\dots,0}_{(n-w)\text{ times}})\in\RR^n$.}$\Tilde{X}\in\RR^n$ such that $\Tilde{X}_i=0$, for all $i>w_{\ell}$.
Again $D$ involves none of the variables $X_i$, where $i\in\{w_{\ell}+1,\dots,n\}$. So each monomial of $DP_2$ involves at least one variable $X_i$, for some $i\in\{w_{\ell}+1,\dots,n\}$, with degree of $X_i$ at least $1$. So $DP_2$ also vanishes at each $\Tilde{X}\in\RR^n$ such that $\Tilde{X}_i=0$, for all $i>w_{\ell}$. 

So for all $(X_1,\dots,X_{w_{\ell}})\in\RR^{w_{\ell}}$, 
since $DP_2(X_1,\dots, X_{w_{\ell}}, 0, \dots, 0) = 0$ we can write
\begin{align*}
    DP_1(X_1,\dots,X_{w_{\ell}})
    &= DP_1(X_1,\dots,X_{w_{\ell}}) + DP_2(X_1,\dots, X_{w_{\ell}}, 0, \dots, 0) \\
    &= DP(X_1,\dots, X_{w_{\ell}}, 0, \dots, 0). 
\end{align*}
The last equality follows from the fact that $P(X_1,\dots,X_n) = P_1(X_1,\dots,X_{w_{\ell}}) + P_2(X_1,\dots,X_n)$.
Next, observe that we can write
\begin{equation*}
\begin{split}
    Q(X_1,\dots,X_{w_{\ell}}) : &= P(1-X_1,\dots, 1-X_{w_{\ell}},0,\dots,0)\\
    &=P_1(1-X_1,\dots,1-X_{w_{\ell}})+ P_2(1-X_1,\dots,1-X_{w_{\ell}},0,\dots,0)\\
    &=P_1(1-X_1,\dots,1-X_{w_{\ell}}).
\end{split}
\end{equation*}
The last equality follows from the fact that $P_2(1-X_1,\dots,1-X_{w_{\ell}},0,\dots,0)=0$. Therefore, for all $(X_1,\dots,X_{w_{\ell}})\in\RR^{w_{\ell}}$, we have
$$
    DQ(X_1,\dots,X_{w_{\ell}})=DP_1(1-X_1,\dots,1-X_{w_{\ell}})
    =DP(1-X_1,\dots,1-X_{w_{\ell}},0,\dots,0).
$$
Now since $P$ vanishes at $\Tilde{a}$ with multiplicity $\ell$, we know that $DP(\Tilde{a})=0$ and hence $DQ(a)=DP(\Tilde{a})=0$. 
%
%
%
Therefore, using Theorem~\ref{prop: Yuval, t,t-1}, we get $\deg(Q)\geq w_{\ell}+2\ell-3$.
\end{proof}

The following proposition shows that the lower bound in Theorem~\ref{th: mult_R} is tight.

 \color{black}

\begin{proposition}
\label{proposition:lower_bound_tightness}
Let $\ell \geq 2$ and $w \in \mathbb{N}$ with $w \geq 2\ell - 3$. Then there exists a polynomial $P \in \mathbb{R}[X_1, \dots, X_n]$ of degree $w + 2\ell - 3$ such that:
\begin{enumerate}
    \item[(a)] 
        $P$ vanishes with multiplicity at least $\ell$ at every $v \in \{0,1\}^n$ with $\w(v) < w$, and

    \item[(b)] 
        For each $\widetilde{w} \in \{w, \dots, n\}$, there exists a $u \in \{0,1\}^n$ with $\w(u) = \widetilde{w}$ such that $P(u) \neq 0$.
\end{enumerate}
\end{proposition}
\color{black}

\begin{proof}
By Theorem~\ref{prop: Yuval, t,t-1}, we know that there exists a polynomial $Q\in\RR[X_1,\dots, X_w]$ of degree $w+2\ell-3$ such that for any $v\in\zone^w\setminus (1,\dots,1)$, $Q$ vanishes at $v$ with multiplicity at least $\ell$ and $Q(1,\dots,1)\neq 0$.

 Now define $P(X_1,\dots,X_n)=Q(X_1,\dots,X_w)$. Then for all $v\in\zone^n$ with 
 $\w( v ) < w$, $P$ vanishes at $v$ and for each $\Tilde{w}\in\{w,\dots,n\}$, the vector $u_{\Tilde{w}}\in\zone^n$ defined by
 $u_{\Tilde{w},i}=1$ if and only if $i\in [\Tilde{w}]$ satisfies $P(u_{\Tilde{w}})=Q(\underbrace{1,\dots,1}_{w\text{ times}})\neq 0$.

 Next consider any differential operator $D=\frac{\partial^{d_1}}{\partial X_1^{d_1}}\dots\frac{\partial^{d_n}}{\partial X_n^{d_n}}$ such that $d_i\geq 0$, for all $i\in [n]$ and $\sum_{i=1}^n d_i< \ell$. We shall show that $DP(v)=0$, for all $v\in\zone^n$ with $\w( v )<w$. First observe that, if $d_i>0$, for some $i> w$, then $DP$ is identically zero. So now we assume that $d_i=0$, for all $i>w$ and take any $v\in\zone^n$ with $\w( v )< w$. Then $\Tilde{v}=(v_1,\dots,v_w)\in\zone^w\setminus\{(1,\dots,1)\}$ and so $\Tilde{v}$ is a zero of $Q$ with multiplicity at least $\ell$. Hence $DQ(\Tilde{v})=0$ and so $DP(v)=DQ(\Tilde{v})=0$, as required.  
\end{proof}

\section*{Acknowledgements}

Arijit Ghosh acknowledges partial support from the Science and Engineering Research Board (SERB), Government of India, through the MATRICS grant MTR/2023/001527, and from the Department of Science and Technology (DST), Government of India, through grant TPN-104427. 
Chandrima Kayal is supported by French PEPR integrated project EPiQ (ANR-22-PETQ-0007).

\bibliographystyle{alpha}

\bibliography{references}

\end{document}